\documentclass[reqno,11pt]{amsart}
\newtheorem{theo}{Theorem}
\newtheorem{rem}{Remark}
\newtheorem{lem}{Lemma}

\newtheorem{prop}{Proposition}
\newtheorem{cor}{Corollary}
\newtheorem{df}{Definition}

\newcommand\eps\varepsilon
\newcommand\ph\varphi
\newcommand\kap\varkappa

\begin{document}

\title[Peano type theorem]
{Peano type theorem for abstract parabolic equations}

\author[Oleg Zubelevich]{Oleg Zubelevich\\ \\\tt
Department of Differential Equations and Mathematical Physics\\
Peoples Friendship University of  Russia\\
Ordzhonikidze st., 3,  117198,  Moscow, Russia\\
 E-mail:  ozubel@yandex.ru}
\email{ozubel@yandex.ru}
\curraddr{2-nd  Krestovskii Pereulok 12-179, 129110, Moscow, Russia}
\thanks{Partially supported by grants RFBR 02-01-00400.}

\begin{abstract}We consider parabolic problems with  non-Lipschitz
nonlinearity in the different scales of Banach spaces and prove local-in-time existence theorem.
New class of parabolic equations that have analytic solutions is obtained.
\end{abstract}

\maketitle
\numberwithin{equation}{section}
\newtheorem{theorem}{Theorem}[section]
\newtheorem{lemma}[theorem]{Lemma}
\newtheorem{definition}{Definition}[section]

\section{Introduction}
This paper is devoted to quasi-linear parabolic equations with a
non-Lipschitz nonlinearity. In the classical setup a
quasi-linear initial value parabolic problem has the form
\begin{equation}
\label{dygr689}
u_t=f(t,u,\nabla^k u)+A u,\quad u\mid_{t=0}=\hat{u}.\end{equation}
Here $A$ is a linear elliptic operator of order $n$ and the term $\nabla^k u$ symbolizes
the derivatives of $u$ up to order $k$.
Besides this, equation (\ref{dygr689}) must be provided with the boundary conditions.

If the function $\hat{u}$ belongs to a suitable space,
 the mapping $f$ is  Lipschitz in a certain sense and $k<n$  then
 problem (\ref{dygr689}) has a unique local-in-time solution. This simple observation easily
follows from the contracting mapping principle.

We consider the case when the function $f$ is non-Lipschitz.
It is well known that in general situation,
in infinitely dimensional  Banach space, an initial value problem for
  differential equation with
non-Lipschitz right hand side does not have solutions \cite{a9,a35,a14}. Nevertheless, as a
rule, the initial value problem lives not in a single Banach space but
in a scale of Banach spaces and in addition this scale is completely
continuous embedded. Such scales for example are the scale of Sobolev spaces, the
scale of analytic functions. This
observation prompts that to find a solution one should
 study the problem in the whole scale.

Let us note another feature of  equations (\ref{dygr689}).
If we reject the Lipschitz hypothesis on $f$ then we obtain a  class of
systems that have existence theorem even in the case when $k\ge n$.
Such type systems remain parabolic in some certain generalized sense.

This effect takes place not only for parabolic equations. If we consider
the Cauchy-Kowalewski problem in the non-Lipschitz setup \cite{zu}
then there are equations such that the order of derivatives in the right
side is greater than in the left one but the solution exists.

These problems do not belong  to the classical partial differential
equations but to the functional-differential equations and the differential equations
with nonlocal terms.

The main mathematical tool we use is a
 locally
convex space version of the Schauder fixed point theorem and
theory of scales of Banach spaces. Another
approaches to  the abstract parabolic problems in   the Lipschitz setup contain in
\cite{arrieta}, \cite{carv}.

\section{Main theorem}
Consider two scales of Banach spaces
$\{E_s,\|\cdot\|^E_s\}_{s>0}$ and $\{G_s,\|\cdot\|^G_s\}_{s>0}$ such that
$E_s\subseteq G_s$ for all $s>0$. All the embeddings
$E_{s+\delta}\subseteq E_s,\quad \delta>0$ are completely continuous and
\begin{equation}
\label{p00}
\|\cdot\|^E_s\le \|\cdot\|^E_{s+\delta}.\end{equation} The parameter $s$
may not  necessarily be  ran through all the positive real numbers.  We do
not use the spaces $E_s,G_s$ with big $s$ and one can assume for example that $s\in (0,1)$.
It is just for simplicity's sake that we consider $s>0.$

Introduce constants
$C,T,R>0,\quad\phi,\alpha\ge 0.$

Let $S^t:G_s\to E_s,\quad t>0$ be a strongly continuous linear semigroup in the following
sense. For any $u\in E_s$ one has
$$\|S^tu-u\|^E_s\to 0\quad \mathrm{as}\quad t\searrow 0
\quad \mathrm{and}\quad \|S^tu\|^E_s\le C\|u\|^E_s.$$
\begin{df}
The semigroup $S^t$ is said to be parabolic if
there exists a constant $\gamma>1$ such that
 for any $\delta,t>0,\quad \delta^\gamma<t<T$
 we have \begin{equation}\label{bvdghr67}\|S^tu\|^E_{s+\delta}\le \frac{C}{t^\phi} \|u\|^G_s.\end{equation}\end{df}

Let $B_s(r)$ be an open ball of the space $E_s$ with radius $r$ and  center
at the origin.
Suppose  a function
$f:(0,T]\times \overline{B}_{s+\delta}(R)\to G_s$ to be  continuous
and such that if $(s+\delta)^\gamma<t\le T$ and $u\in\overline{B}_{s+\delta}(R)$
  then the following inequality holds
\begin{equation}
\label{fuck_you_putin}
\|f(t,u)\|^G_s\le
\frac{C}{\delta^\alpha}.\end{equation}
\begin{rem}\label{dtby567}A case when
$$\|f(t,u)\|^G_s\le
\frac{C}{t^\beta\delta^\alpha},\quad \beta>0$$ is rather usual but since
$\delta^\gamma<t$ this case reduces to (\ref{fuck_you_putin}):
$C/(t^\beta\delta^\alpha)\le C/\delta^{\beta\gamma+\alpha}.$\end{rem}

We proceed with two setups of our problem. The first one is a classical setup
and we find classical solutions and the second one is a generalized setup
to obtain generalized solutions.

In the generalized setup we are looking for solutions to the following
integral equation
\begin{equation}\label{ftbyr6}u(t)=\int_0^tS^{(t-\xi)}f(\xi,u(\xi))\,d\xi.
\end{equation}

In the classical setup we make several additional assumptions. Namely,
suppose that $G_s=E_s$. Introduce  a linear operator
 $A:E_{s+\delta}\to E_s$ and assume that the semigroup $S^t$ is generated by this
 operator: $S^t=e^{At}$ such that
 for any $u\in E_{s+\delta}$ we have
\begin{equation}\label{dyh6}
\lim_{h\to
0+}\Big\|\frac{1}{h}\Big(e^{Ah}-\mathrm{id}_{E_{s+\delta}}\Big)u-Au\Big\|^E_s=0.\end{equation}

In the classical setup our problem has the form
\begin{align}
u_t&=f(t,u)+Au,\label{main_eq}\\
 u\mid_{t=0}&=0\label{main_eq11}.
\end{align}
The sense of initial condition (\ref{main_eq11}) will be clear in the sequel.

Now we give a definition.
\begin{df}
We shall say that problem (\ref{main_eq}) or (\ref{ftbyr6}) is parabolic if
the semigroup $S^{t}$ is parabolic and $$\chi=\phi+\frac{\alpha}{\gamma}<1.$$
\end{df}In case of remark \ref{dtby567} $\chi=\phi+\beta+\alpha/\gamma.$

Let a space $E^1(T),\quad T>0$ be given by the formula
\begin{equation}
\label{pewrpwe4}
E^1(T)=\bigcap_{0<s^\gamma<\tau<T}C^1((\tau,T),E_{s}).\end{equation}
This space consists of all functions $u$ that map any number $t\in (0,T)$
to the element $u(t)\in \bigcap_{0<s^\gamma<t}E_s$ and  the
restriction $u\mid_{(\tau,T)}$ belongs to the space $C^1((\tau,T),E_{s})$
for all $s\in (0,\tau^{1/\gamma})$.

\begin{theo}
\label{main_th}
1) Classical setup. Suppose that problem (\ref{main_eq})
is parabolic.
Then there exists a  constant $T_*>0$ such that  this problem has
a solution
$u(t)\in E^1(T_*),$
and for any constant $c\in (0,1)$ one has
\begin{equation}
\label{ab} \|u(t)\|^E_{ct^{1/\gamma}}\to 0\quad\mathrm{as}\quad t\searrow 0.
\end{equation} The function $u(t)$ also solves equation $(\ref{ftbyr6})$.

2) Generalized setup. Suppose that problem (\ref{ftbyr6})
is parabolic.
Then there exists a  constant $T_*>0$ such that  this problem has
a solution
$$u(t)\in E(T_*)=\bigcap_{0<s^\gamma<\tau<T_*}C((\tau,T_*),E_{s}).$$

In both cases the constants $T_*$ depends only on $C,\alpha,\gamma,\phi.$
\end{theo}
The proof of theorem \ref{main_th} contains in sections \ref{fhh234556}, \ref{oiuio}.

Then to illustrate the effect discussed in the
Introduction,
 theorem \ref{main_th} is applied
to  a nonlocal parabolic problem.
To compare our result with the known one we also consider the Navier-Stokes equation.

If
$A$ is the classical Laplace operator and  the parabolic equation is
considered in a suitable domain then $\gamma=2$ and the inequality from formula
(\ref{pewrpwe4}) takes the form $0<s^2<\tau$.

The parameter $s$ symbolizes a spatial variable, so
that this inequality specifies the parabolic domain in the
plane $(\tau,s)$. This endows the term "parabolic equation" with the new
sense.

Let us remark that if $G_s=E_s=\mathbb{R}^m,\quad \|\cdot\|^E_s=|\cdot|,\quad s>0$ and
$A=0$ then theorem \ref{main_th} generalizes classical Peano's theorem to
the case when the right side of the equation satisfies
(\ref{fuck_you_putin}) with
$s=\delta=(t/3)^{1/\gamma}.$

\section{Preliminaries on functional analysis}
\label{fhh234556}
In this section we collect several facts from functional analysis. These
facts will be useful in the section \ref{oiuio} when we prove theorem \ref{main_th}.


Consider the spaces $$C([\tau,T],E_{\mu\tau^{1/\gamma}}),\quad 0<\mu<1,\quad 0<\tau<T$$ with
standard norms. Now we construct the projective limit of these spaces.
 Define a space $E(T)$ as follows
$$E(T)=\bigcap_{0<\mu<1}\bigcap_{0<\tau<T}C([\tau,T],E_{\mu\tau^{1/\gamma}}).$$
There is another equivalent definition of  the space $E(T)$:
$$E(T)=\bigcap_{0<s^\gamma<\tau<T}C([\tau,T],E_{s}).$$

Being endowed with a collection of seminorms
\begin{equation}
\label{hjjk}
\|u\|_{\tau,\mu}=\max_{\tau\le \xi\le
T}\|u(\xi)\|^E_{\mu\tau^{1/\gamma}},\quad u\in E(T)\end{equation}
the space $E(T)$ becomes a locally convex topological space.

  These seminorms obviously satisfy the following inequalities
\begin{align}
\|u\|_{\tau,\mu}&\le \|u\|_{\tau,\mu+\delta},\quad \delta>0,\label{o1}\\
\|u\|_{\tau,r\mu}&\le \|u\|_{r^\gamma\tau,\mu},\quad 0<r\le 1.\label{o2}\end{align}
Indeed, formula (\ref{o1}) follows from (\ref{p00}) directly.
Formula (\ref{o2}) is a result of the estimate
$$\|u\|_{\tau,r\mu}=\max_{\tau\le \xi\le
T}\|u(\xi)\|^E_{\mu(r^\gamma\tau)^{1/\gamma}}\le\max_{r^\gamma\tau\le \xi\le
T}\|u(\xi)\|^E_{\mu(r^\gamma\tau)^{1/\gamma}}=\|u\|_{r^\gamma\tau,\mu}.$$

Formulas (\ref{o1}), (\ref{o2}) imply that the space $E(T)$ is first
countable: the   topology of this space can be defined  by the seminorms
(\ref{hjjk}) only with $\mu,\tau\in \mathbb{Q}$.

Recall the  Arzela-Ascoli theorem \cite{Schwartz}:
\begin{theo}
\label{Ascoli}
Let $H\subset C([0,T],X)$ be
a set in the space of continuous functions
with values in a Banach space $X$.
Assume that the set $H$ is closed, bounded, uniformly continuous and
for every $t\in [0,T]$ the set $\{u(t)\in X\}$ is a compact set in the space
$X$. Then the set $H$ is a compact set in the space $C([0,T],X)$.
\end{theo}
Now we shall establish an analogue of this result.
\begin{prop}
\label{yun}
Suppose that a set $K\subset E(T)$ is closed. Then $K$ is a compact set if the
following two conditions are fulfilled.

The set $K$ is bounded.

For any $\eps>0$ and for any $\tau\in (0,T),\quad \mu\in(0,1)$ there
exists a constant $\delta>0$ such that if $t',t''\in [\tau,T],\quad
|t'-t''|<\delta$ then
$$\sup_{u\in K}\|u(t')-u(t'')\|^E_{\mu\tau^{1/\gamma}}<\eps.$$
(This means that $K$ is a uniformly continuous set.)
\end{prop}
First prove a lemma.
\begin{lem}
\label{vvdf}
Let $\{v_j\}\subseteq K$ be a sequence. Then for any  $\tau\in (0,T)$ the sequence
$\{v_j\}$ contains a subsequence that is convergent in all the norms
$\|\cdot\|_{\tau,\mu},\quad\mu\in (0,1).$\end{lem}
\proof
Indeed, take an increasing sequence $\mu_k\to 1,\quad \mu_1>0$ and fix any
value of $\tau\in (0,T)$. Since the sequence $\{v_j\}$
is bounded and uniformly continuous in
$C([\tau,T],E_{\mu_2\tau^{1/\gamma}})$ then by theorem \ref{Ascoli} it  contains
a subsequence
$\{v^1_j\}$ that is convergent
 in $C([\tau, T],E_{\mu_1\tau^{1/\gamma}})$.

 Further since the sequence $\{v^1_j\}$ is bounded and
uniformly continuous in $C([\tau, T],E_{\mu_3\tau^{1/\gamma}})$ one can
pick  a subsequence $\{v^2_j\}\subseteq \{v^1_j\}$ such that the
sequence $\{v^2_j\}$ is convergent in
$C([\tau, T],E_{\mu_2\tau^{1/\gamma}})$ etc.

 By inequality (\ref{o1})  the diagonal
sequence $\{v^j_j\}$ converges in all the norms
$\|\cdot\|_{\tau,\mu},\quad \mu\in (0,1)$ with this fixed $\tau$.\endproof

{\it Proof of proposition \ref{yun}.}
A set $P=\mathbb{Q}\bigcap (0,T)$ is countable. So we can number its
elements as follows $P=\{\tau_i\}_{i\in \mathbb{N}}$.

We must show that any sequence $\{u_j\}\subseteq K$ contains a convergent
subsequence $\{u_{j_k}\}.$

By lemma \ref{vvdf} there is a subsequence $\{u^1_j\}\subseteq\{u_j\}$
that is convergent in all the norms $\|\cdot\|_{\tau_1,\mu}\quad \mu\in (0,1).$
By the same argument there is a subsequence $\{u^2_j\}\subseteq \{u^1_j\}$
that is convergent in all the norms $\|\cdot\|_{\tau_2,\mu}\quad \mu\in
(0,1)$ etc.

The diagonal sequence $\{u^j_j\}$ is convergent in all the
norms $\|\cdot\|_{\tau_k,\mu},\quad k\in \mathbb{N},\quad \mu\in (0,1).$

By inequality (\ref{o2}) the sequence $\{u^j_j\}$ is convergent in all the
norms $\|\cdot\|_{\tau,\mu},\quad \tau\in (0,T),\quad \mu\in (0,1).$

Proposition \ref{yun} is proved.

\begin{lem}
\label{gt34556}
Let $X,Y$ be Banach spaces. Suppose that $A_a:X\to Y,\quad a'>a>0$ is a collection
of bounded
linear operators such that for each $x\in X$ we have
$$\sup_{a'>a>0}\|A_a x\|_{Y}<\infty,\quad
\|A_ax\|_{Y}\to 0\quad\mathrm{as}\quad a\to 0.$$

 Then for any compact set
$B\subset X$ it follows that
$$\sup_{x\in B}\|A_ax\|_{Y}\to 0\quad \mathrm{as}\quad a\to 0.$$
\end{lem}This result is a direct consequence of the Banach-Steinhaus theorem \cite{Schwartz}.

Let us recall  a generalized version of the Schauder fixed
point theorem.
\begin{theo}[\cite{Browder}]\label{Schau-Tich}Let $W$ be a closed convex subset of the
locally convex space $E$. Then a compact continuous mapping $f:W\to W$ has a fixed
point $\hat u$ i.e. $f(\hat u)=\hat u$.\end{theo}

\section{Proof of Theorem \protect\ref{main_th}}\label{oiuio}

By definition put
$$W(T_*)=\{u\in E(T_*)\mid \|u\|_{\tau,\nu}\le R,\quad 0<\tau<T_*,\quad 0<\nu<1\}.$$
The constant $T_*>0$ will be defined.

First we  find a fixed point of a mapping
$$F(u)=\int_0^tS^{t-\xi}f(\xi,u(\xi))\,d\xi.$$ This fixed point is
 the generalized solution announced in the second part of the theorem.
Then by using formula (\ref{dyh6}) we show that this fixed point is the desired solution to problem
(\ref{main_eq}).

\begin{lem}
\label{lk90}If the constat $T_*$ is small enough then the mapping $F$ takes the
set $W(T_*)$ to itself.\end{lem}
\proof Let  constants $t,s$ be taken as follows $0<s< t^{1/\gamma},\quad t\le T_*.$
Suppose $u\in W(T_*)$ then estimate a function $v(t)=F(u):$
\begin{equation}\label{gfd6}
\|v(t)\|^E_{s}\le\int_0^t \|S^{t-\xi}f(\xi,u(\xi))\|^E_s\,d\xi=X+Y,
\end{equation}
here we use the notation
$$
X=\int_0^{t-s^\gamma} \|S^{t-\xi}f(\xi,u(\xi))\|^E_s\,d\xi,\quad
Y=\int_{t-s^\gamma}^t
\|S^{t-\xi}f(\xi,u(\xi))\|^E_s\,d\xi.$$
To estimate $X$ take   constants $\eps$ and  $\mu$ such that
\begin{equation}
\label{poo0999}
0<\eps<\frac{s}{t^{1/\gamma}}<\mu<1.\end{equation} The constant $\eps$ is assumed to be small and
the constant $\mu$ is assumed to be close to $1$.

Let the  variables $\delta$ and $\delta'$ be given by the formulas
$$\delta=s-\eps\xi^{1/\gamma},\quad \delta'=\xi^{1/\gamma}(\mu-\eps).$$
Taking into account that $\xi \in (0,t-s^\gamma]$ we  see that the
variables $\delta ,\delta'$ are positive and
\begin{equation}
\label{gt6}
s-\delta>0,\quad s-\delta+\delta'<\xi^{1/\gamma},\quad
\delta<(t-\xi)^{1/\gamma}.\end{equation}
The inequality in the middle implies that
\begin{equation}
\label{09ii8u}
u(\xi)\in \overline{B}_{s-\delta+\delta'}(R)\end{equation}
and thus the term $X$ is estimated as follows
\begin{align}
X&\le C\int_0^{t-s^\gamma}(t-\xi)^{-\phi}\|f(\xi,u(\xi))\|^G_{s-\delta}\,d\xi
\le
C^2\int_0^{t-s^\gamma}\frac{1}{\delta'^\alpha(t-\xi)^{\phi}}
\,d\xi\nonumber\\
&\le \frac{C^2}{(\mu-\eps)^\alpha}
\int_0^{t-s^\gamma}\frac{d\xi}{(t-\xi)^{\phi}\xi^{\alpha/\gamma}}\Big|_{\xi=yt}=
\frac{C^2t^{1-\chi}}{(\mu-\eps)^\alpha}\int_0^{1-s^\gamma/t}\frac{dy}{(1-y)^\phi
y^{\alpha/\gamma}}\nonumber\\
&\le \frac{C^2Jt^{1-\chi}}{(\mu-\eps)^\alpha},\quad J=\int_0^{1}\frac{dy}{(1-y)^\phi
y^{\alpha/\gamma}}.
\label{tyu6}\end{align}
We shall estimate the term $Y$.

Introduce a function $\psi$ by the formula
$$\psi(y)=y^{1/\gamma}+(1-y)^{1/\gamma}-1.$$ The function $\psi$ is
positive on the interval $(0,1)$.
Define a constant $I$ as follows
$$I=\int_0^1\frac{dy}{(1-y)^{\phi}(\psi(y))^\alpha}.$$
Let the constant $\mu$ be as above. We redefine the variables
$\delta,\delta'$ by the formulas
$$\delta=\mu(t-\xi)^{1/\gamma},\quad \delta'=\mu\xi^{1/\gamma}+\delta-s.$$
Now the variable $\xi$ belongs to the interval $[t-s^\gamma,t]$ and thus the
variables $\delta,\delta'$ are positive and satisfy inequalities
(\ref{gt6}).

It is only not trivial to show that the variable $\delta'$ is
positive. Let us prove this. Indeed,
\begin{equation}
\label{rty78}
\delta'=\mu\xi^{1/\gamma}+\mu(t-\xi)^{1/\gamma}-s
=t^{1/\gamma}\Big(\mu
y^{1/\gamma}+\mu(1-y)^{1/\gamma}-\frac{s}{t^{1/\gamma}}\Big),\end{equation}
recall that
$y=\xi/t.$
Form (\ref{rty78}) it follows that
\begin{equation}
\label{poo099}
\delta'>t^{1/\gamma}\mu\psi(y).\end{equation}
By the same argument as above, inclusion
 (\ref{09ii8u}) is fulfilled with the new $\delta$ and $\delta'$.

We are ready to estimate the term $Y$. By (\ref{poo099})
it follows that
\begin{align}
Y&\le C\int_{t-s^\gamma}^t(t-\xi)^{-\phi}\|f(\xi,u(\xi))\|^G_{s-\delta}\,d\xi \le C^2
\int_{t-s^\gamma}^t \frac{d\xi}{(t-\xi)^{\phi}\delta'^\alpha}\nonumber\\
&\le \frac{C^2t^{1-\chi}}{\mu^\alpha}
\int_{1-s^\gamma/t}^1\frac{dy}{(1-y)^{\phi}(\psi(y))^\alpha}\le
\frac{C^2I}{\mu^\alpha}t^{1-\chi}.\label{fghfy879879}
\end{align}
Now the assertion the of lemma follows from formulas (\ref{gfd6}), (\ref{tyu6}) and
(\ref{fghfy879879}).
\endproof
\begin{cor}\label{fgt674}Formulas (\ref{tyu6}),
(\ref{fghfy879879}) imply that if $0<s^\gamma<t\le T_*$ and $v(t)=F(u),\quad u\in W(T_*)$ then
$$
\|v(t)\|^E_s\le c_2 t^{1-\chi},$$
here $c_2$ is a positive constant independent on $u,t,s.$\end{cor}

\begin{lem}
\label{ytyt67}
The set $F(W(T_*))$ is precompact in $E(T_*)$.
\end{lem}
\proof By proposition \ref{yun} it is sufficient to prove that the set
$F(W(T_*))$ is uniformly continuous.

Take a function $u\in W(T_*)$ and let $v(t)=F(u)$.
We must show that if $t',t''\ge\tau,\quad \tau\in (0,T_*)$ then for any $\mu\in (0,1)$ one has
$$\sup_{u\in W(T_*)}\|v(t')-v(t'')\|^E_{\mu\tau^{1/\gamma}}\to 0,\quad \mathrm{as}\quad
|t'-t''|\to 0.$$

Indeed, for definiteness assume that $t''>t'$ then
\begin{align}
v(t'')-v(t')&=\int_{t'}^{t''}S^{t''-\xi}f(\xi,u)\,d\xi\nonumber\\
&+\Big(S^{t''-t'}-\mathrm{id}_{E_s}\Big)
\int_0^{t'}S^{t'-\xi}f(\xi,u)\,d\xi,\quad s^\gamma<\tau.\label{nm75}\end{align}
Choose a positive constant $\delta$ such that $(s+\delta)^\gamma<\tau$
and using the parabolicity of the semigroup $S^t$
  estimate the first term from the right side of this formula
\begin{align}
\Big\|\int_{t'}^{t''}S^{t''-\xi}f(\xi,u)\,d\xi\Big\|^E_s&\le
C \int_{t'}^{t''}(t''-\xi)^{-\phi}\|f(\xi,u)\|^G_s\,d\xi
\nonumber\\&\le
C^2 \int_{t'}^{t''}\frac{d\xi}{\delta^\alpha(t''-\xi)^\phi}
=\frac{C^2}{\delta^\alpha(1-\phi)}(t''-t')^{1-\phi}.\nonumber
\end{align}
 So that the first
term in the right side of (\ref{nm75}) is vanished uniformly.

Consider a set
$$U=\bigcup_{\tau\le t'\le T_*}\Big\{\int_0^{t'}S^{t'-\xi}f(\xi,u)\,d\xi\mid u\in
W(T_*)\Big\}.$$  By lemma \ref{lk90} the set $U$ is bounded in any space
$E_{\mu'\tau^{1/\gamma}}$ with $1>\mu'>\mu$ thus it is compact in
$E_{\mu\tau^{1/\gamma}}$.
  By lemma \ref{gt34556}
we get
$$\sup_{w\in U}\|S^{t''-t'}w-w\|^E_{\mu\tau^{1/\gamma}}\to
0,\quad \mathrm{as}\quad t''-t'\to 0.$$
This shows that the second term in the right side of formula (\ref{nm75}) is vanished uniformly.
\endproof
\begin{cor}\label{cvb41}
The set $F(W(T_*))$ is uniformly continuous with respect to the variable
$t$.\end{cor}

\begin{lem}\label{vczq3}The mapping $F:W(T_*)\to W(T_*)$ is continuous with respect to
the topology of the space $E(T_*)$.\end{lem}
\proof
Suppose a sequence $\{v_l\}\subset W(T_*)$ to be convergent to the element
$v\in W(T_*)$ as $l\to \infty.$
We need to show that for any $s^\gamma<\tau<T_*$ the sequence
$$\sup_{\tau\le t\le
T_*}\Big\|\int_0^tS^{t-\xi}f(\xi,v_l(\xi))\,d\xi-
\int_0^tS^{t-\xi}f(\xi,v(\xi))\,d\xi\Big\|^E_s$$
vanishes as $l\to \infty.$

By corollary \ref{cvb41} the sequence
\begin{equation}
\label{abcd3453}
\Big\{\int_0^tS^{t-\xi}f(\xi,v_l(\xi))\,d\xi\Big\}\end{equation}
 is uniformly continuous on the interval $[\tau,T_*]$.
  The uniform convergence of such a sequence is equivalent to
its pointwise convergence \cite{Schwartz}. Thus it is sufficient to prove that sequence
(\ref{abcd3453}) is convergent in $E_s$ for each $t\in [\tau,T_*]$.

Fix $t\in [\tau,T_*]$ and let  constants $\eps,\mu$ satisfy inequality (\ref{poo0999}).
Then using the argument of lemma \ref{lk90} write
\begin{align}
\Big\|&\int_0^tS^{t-\xi}(f(\xi,v_l(\xi))-f(\xi,v(\xi)))\,d\xi\Big\|^E_s\nonumber\\
&\le
\int_0^{t-s^\gamma}(t-\xi)^{-\phi}\|f(\xi,v_l(\xi))-f(\xi,v(\xi))\|^G_{\eps\xi^{1/\gamma}}\,d\xi\nonumber\\&+
\int_{t-s^\gamma}^t(t-\xi)^{-\phi}\|f(\xi,v_l(\xi))-f(\xi,v(\xi))\|^G_{s-\mu(t-\xi)^{1/\gamma}}\,d\xi.\label{vbcbvc}
\end{align}
Since the function $f$ is continuous, for a fixed $\xi$ we have:
\begin{align}(t-\xi)^{-\phi}\|f(\xi,v_l(\xi))-f(\xi,v(\xi))\|^G_{\eps\xi^{1/\gamma}}&\to 0,
\quad\xi\in[0,t-s^\gamma],\nonumber\\
(t-\xi)^{-\phi}\|f(\xi,v_l(\xi))-f(\xi,v(\xi))\|^G_{s-\mu(t-\xi)^{1/\gamma}}&\to 0,\quad \xi\in[t-s^\gamma,t),
\nonumber\end{align}as $ l\to
\infty.$

Moreover by formulas (\ref{tyu6}), (\ref{fghfy879879})
both of these expressions are majorized with the $L^1$-integrable
function:
\begin{align}
(t&-\xi)^{-\phi}\|f(\xi,v_l(\xi))-f(\xi,v(\xi))\|^G_{\eps\xi^{1/\gamma}}\nonumber\\&\le
(t-\xi)^{-\phi}(\|f(\xi,v_l(\xi))\|^G_{\eps\xi^{1/\gamma}}+\|f(\xi,v(\xi))\|^G_{\eps\xi^{1/\gamma}})
\le
\frac{2C^2}{(\mu-\eps)^\alpha\xi^{\alpha/\gamma}(t-\xi)^{\phi}}\nonumber,
\end{align}
and
$$(t-\xi)^{-\phi}\|f(\xi,v_l(\xi))-f(\xi,v(\xi))\|^G_{s-\mu(t-\xi)^{1/\gamma}}\le
\frac{2C^2}{ t^{\alpha/\gamma}\mu^\alpha(\psi(\xi/t))^\alpha(t-\xi)^{\phi}}.$$
Therefore by the Dominated convergence theorem the integrals in the right
side of (\ref{vbcbvc}) are vanished as $l\to \infty$.
\endproof

So by theorem \ref{Schau-Tich} and  lemmas \ref{lk90}, \ref{ytyt67}, \ref{vczq3}
we obtain a fixed point of the mapping $F$, say $u$:
$$F(u)=u\in W(T_*).$$
This proves the second part of theorem \ref{main_th}.


To prove the first one
let us show that this fixed point  is the  solution to
problem (\ref{main_eq}).
Suppose that $t,t+h>s^\gamma$. First consider the case $h>0$.
 Differentiate the function $u(t)$ explicitly:
\begin{align}
u_t(t)&=\lim_{h\to
0}h^{-1}\Big(\int_0^{t+h}e^{A(t+h-\xi)}f(\xi,u(\xi))\,d\xi-
\int_0^te^{A(t-\xi)}f(\xi,u(\xi))\,d\xi\Big)\nonumber\\
&=\lim_{h\to
0}h^{-1}\int_t^{t+h}e^{A(t+h-\xi)}f(\xi,u(\xi))\,d\xi\nonumber\\
&+
\lim_{h\to
0}h^{-1}(e^{Ah}-\mathrm{id}_{E_s})\int_0^te^{A(t-\xi)}f(\xi,u(\xi))\,d\xi.
\label{njyhgu645}\end{align}

Lemma \ref{lk90} implies that
$\int_0^te^{A(t-\xi)}f(\xi,u(\xi))\,d\xi\in E_{s'}$ with $s^\gamma<s'^\gamma<t,t+h$ hence
  formula (\ref{dyh6}) gives
\begin{equation}
\label{ps2}
h^{-1}(e^{Ah}-\mathrm{id}_{E_s})\int_0^te^{A(t-\xi)}f(\xi,u(\xi))\,d\xi\to
A\int_0^te^{A(t-\xi)}f(\xi,u(\xi))\,d\xi\end{equation} in $E_{s}$  as $h\to
0.$

Let us prove that
\begin{equation}\label{ps1}
h^{-1}\int_t^{t+h}e^{A(t+h-\xi)}f(\xi,u(\xi))\,d\xi\to
f(t,u(t))\end{equation} in $E_s$ as $h\to 0$.

Indeed, observe that
\begin{align}
h^{-1}&\int_t^{t+h}e^{A(t+h-\xi)}f(\xi,u(\xi))\,d\xi-f(t,u(t))\nonumber\\&=
h^{-1}\Big(\int_t^{t+h}e^{A(t+h-\xi)}(f(\xi,u(\xi))-f(t,u(t)))\,d\xi\nonumber\\&+
\int_t^{t+h}(e^{A(t+h-\xi)}-\mathrm{id}_{E_s})f(t,u(t))\,d\xi\Big).\nonumber\end{align}
The first integral in the right side of this formula is estimated as
follows:
\begin{align}
\Big\|\int_t^{t+h}&e^{A(t+h-\xi)}(f(\xi,u(\xi))-f(t,u(t)))\,d\xi\Big\|^E_s\nonumber\\&\le
Ch\max_{t\le\xi\le t+h}\|f(\xi,u(\xi))-f(t,u(t))\|^E_s=o(h).\nonumber
\end{align}
Since the semigroup $e^{At}$ is strongly continuous for the second integral we
get
\begin{align}
\Big\|\int_t^{t+h}&(e^{A(t+h-\xi)}-\mathrm{id}_{E_s})f(t,u(t))\,d\xi\Big\|^E_s\nonumber\\&\le
h\max_{t\le\xi\le t+h}\|(e^{A(t+h-\xi)}-\mathrm{id}_{E_s})f(t,u(t))\|^E_s=o(h).\nonumber
\end{align}

If $h<0$  then instead of formula (\ref{njyhgu645}) one must use the
following expression
\begin{align}
u_t(t)&=\lim_{h\to
0}h^{-1}\Big((\mathrm{id}_{E_s}-e^{-Ah})\int^{t+h}_{0}e^{A(t+h-\xi)}f(\xi,u(\xi))\,d\xi\nonumber\\
&-\int_{t+h}^te^{A(t-\xi)}
f(\xi,u(\xi))\,d\xi\Big)\nonumber.
\end{align}
In this case only the proof of the formula
\begin{align}\lim_{h\to
0}h^{-1}(\mathrm{id}_{E_s}&-e^{-Ah})\int^{t+h}_{0}e^{A(t+h-\xi)}f(\xi,u(\xi))\,d\xi\nonumber\\
&=A\int^{t}_{0}e^{(t-\xi)}f(\xi,u(\xi))\,d\xi\nonumber\end{align}
differs from the previous argument.

Let us prove this formula.
Obviously we have
\begin{align}
(\mathrm{id}_{E_s}&-e^{-Ah})\int^{t+h}_{0}e^{A(t+h-\xi)}f(\xi,u(\xi))\,d\xi\nonumber\\&=
(\mathrm{id}_{E_s}-e^{-Ah})u(t)+(\mathrm{id}_{E_s}-e^{-Ah})(u(t+h)-u(t)).\label{ppovb7}\end{align}
The set $$V=\Big\{\frac{u(t+h)-u(t)}{\|u(t+h)-u(t)\|^E_{s'}}\Big| h\in (h',0)\Big\}$$
with $h'<0$ close to zero  is bounded in $E_{s'},\quad s^\gamma<s'^\gamma<t+h'$.
Consequently
$V$ is a compact set in $E_s$. By lemma \ref{gt34556} the set
$$(A_{-h}-A)V,\quad
A_{-h}=\frac{1}{h}\Big(\mathrm{id}_{E_s}-e^{-Ah}\Big)$$ is bounded in $E_s$
and thus the set $A_{-h}V$ is also bounded.

Thus taking into account that the
  function $u(t)$ is continuous we yield
\begin{align}\Big\|\frac{1}{h}\Big(\mathrm{id}_{E_s}&-e^{-Ah}\Big)(u(t+h)-u(t))\Big\|^E_s\nonumber\\
 &=
\|u(t+h)-u(t)\|^E_{s'}\cdot\Big\|A_{-h}
\frac{u(t+h)-u(t)}{\|u(t+h)-u(t)\|^E_{s'}}\Big\|^E_{s}
=o(1).\nonumber\end{align}
For the second term of the right side of (\ref{ppovb7}) this implies
$$\|(\mathrm{id}_{E_s}-e^{-Ah})(u(t+h)-u(t))\|^E_s=o(h).$$
The first term of the right side of formula (\ref{ppovb7}) is estimated as
follows
$$\|(\mathrm{id}_{E_s}-e^{-Ah})u(t)-hAu(t)\|^E_s=o(h).
$$

Substituting formulas (\ref{ps2})  and (\ref{ps1}) to (\ref{njyhgu645}) we
see that the function $u$ is a solution to equation (\ref{main_eq}).

Formula (\ref{ab}) follows from corollary \ref{fgt674}.

Theorem \ref{main_th} is proved.

\section{Applications}In the sequel we denote  all the
inessential positive constants by the same letter
$c$.

\subsection{Parabolic equation with gradient nonlinearity}In
this section we consider a model example.

Let $M\subset \mathbb{R}^m$  be a bounded domain with smooth boundary $\partial M$.

Consider the following equation
\begin{equation}
\label{heateq1}
u_t=f(\nabla u)+\Delta u,\quad u\mid_{t=0}=\hat u\in H_0^{1,q}(M),\quad u(t,\partial M)=0,\quad
t>0,\end{equation} here $q>1$.

The function $f$ is
continuous in
$\mathbb{R}^m$ and for
all $z\in \mathbb{R}^m$ we have $|f(z)|\le c(|z|^p+1),\quad  q\ge p\ge 1$.
Note that the function $f$  may not
necessarily be a Lipschitz function.

Let us show that if
\begin{equation}
\label{subcritfgh}
m(p-1)<q.\end{equation}
then
problem (\ref{heateq1}) has a generalized
solution from $C([0,T],H_0^{1,q}(M)),$ the constant $T>0$
depends on $\hat u$.

If the function  $f$ is a Lipschitz function then inequality (\ref{subcritfgh})
is well known: it  corresponds to the subcritical case in the
sense of Fujita.

After the change of the unknown function $u=e^{\Delta t}\hat{u}+v$ our
problem takes the form
\begin{equation}
\label{heateq2}
v_t=g(t,x,\nabla v)+\Delta v,\quad v\mid_{t=0}=0,\quad g(t,x,\nabla v)=
f(\nabla(e^{\Delta t}\hat{u}+v)).\end{equation}

Consider  problem (\ref{heateq2}) in the scales
$$E_s=H_0^{1+s_0+s,q}(M),\quad \|\cdot\|_s=\|\cdot\|_{H_0^{1+s_0+s,q}(M)},\quad
s\in (0,S),$$
 and $$G_s=H^{-\lambda,q}(M),\quad \|\cdot\|^G=\|\cdot\|_{H^{-\lambda,q}(M)},$$
 this means that all the spaces $G_s$ coincide with each other, the constants $S>0,s_0\ge
0$ and $0\le\lambda<m(1-1/q)$ to be defined.

Introduce a constant
$$
r=\frac{qm}{m+\lambda q}\in (1,q].$$
Then using standard facts on the Sobolev spaces estimate the function $g$:
\begin{align}\|g(t,x,\nabla v)\|^G&\le c\|g(t,x,\nabla
v)\|_{L^r(M)}\le c(\|\nabla(e^{\Delta t}\hat{u}+v)\|^p_{L^{pr}(M)}+1)\nonumber\\
&\le c(\|e^{\Delta t}\hat{u}\|^p_{H^{1,pr}(M)}+\|v\|^p_{H^{1,pr}(M)}+1).\label{dgy678}\end{align}
Choose a constant $s_0$ as follows
$$s_0=m\Big(\frac{1}{q}-\frac{1}{rp}\Big).$$
 Then the condition $H_0^{1+s_0,q}(M)\subseteq
H_0^{1,pr}(M)$ is satisfied.

Here we assume that the constant $\lambda$ is such that we have $q<rp$.
Note that
$$\|e^{\Delta t}\hat{u}\|^p_{H^{1,pr}(M)}\le
ct^{-\beta}\|\hat{u}\|^p_{H^{1,q}(M)},\quad\beta=\frac{m}{2}\Big(\frac{p}{q}-\frac{1}{r}\Big).$$
If $v\in B_s=\{h\in E_s\mid\|h\|^E_s\le 1\}$ then by all these argument formula  (\ref{dgy678})
implies
$$\|g(t,x,\nabla v)\|^G
\le \frac {c}{t^\beta}.$$

Another inequality we need is
$$\|e^{\Delta t}w\|_s\le c t^{-\phi}\|w\|^G,\quad \phi=\frac{1+s_0+s+\lambda}{2},$$
this formula also follows from the standard Sobolev spaces theory.

\begin{prop}
The mapping $(t,v)\mapsto g(t,x,\nabla v)$ is a continuous mapping of
$(0,T)\times B_s$ to $G_s$.\end{prop}

\proof Assume the converse: there exists a sequence
$(t_k,v_k)$ such that $t_k\to t\in (0,T),\quad v_k\to v$ in $E_s$ as
$k\to\infty,\quad v,v_k\in B_s$ and
\begin{equation}
\label{ft7yhn}\|g_k(x)-g(x)\|^G\ge c>0,\end{equation}
here we put $g_k(x)=g(t_k,x,\nabla v_k),\quad g(x)=g(t,x,\nabla v)$.

By the argument above formula (\ref{ft7yhn}) imply
$$\|g_k(x)-g(x)\|_{L^r(M)}\ge c>0,$$
Since  $\nabla v_k\to \nabla v$ in $L^{pr}(M)$ then
there exists a subsequence $\{v_{k'}\}\subseteq  \{v_{k}\}$ such that
$\nabla v_{k'}\to \nabla v$ almost every where in $M$.
Thus $|g_{k'}(x)-g(x)|^r\to 0$ almost everywhere in $M$. Consequently
$|g_{k'}(x)-g(x)|^r\to 0$ in measure.

It remains to show that the sequence $|g_{k'}(x)-g(x)|^r$ is uniformly
integrable. If we do this then by the Vitali convergence theorem \cite{fol} it
follows that
$\|g_{k'}(x)-g(x)\|_{L^r(M)}\to 0$ and this contradiction proves the Proposition.

Note that since $v_k,v \in E_s\quad s>0$ we actually have
  $\nabla v_k\to \nabla v$, in $L^{pr+\sigma}(M)$ with small $\sigma>0$.
 Thus
 the functions $g_{k'}(x)-g(x)$ belong not  only
to $L^r(M)$ but also to $L^{r+\eps}(M)$ with  small $\eps>0$ and
the sequence
$\|g_{k'}(x)-g(x)\|_{L^{r+\eps}(M)}$ is bounded (these observations follow
 from the same argument as above). The last observation  can
be rewritten as follows:
\begin{align}\sup_{k'}&\int_M|g_{k'}(x)-g(x)|^r\ae(|g_{k'}(x)-g(x)|)\,dx\nonumber\\
&=\sup_{k'}\|g_{k'}(x)-g(x)\|^{r+\eps}_{L^{r+\eps}(M)}
<\infty\nonumber\end{align}
with $\ae(y)=y^\eps$. Since the function $\ae$ is monotone and unbounded in $\mathbb{R}_+$, this proves
the uniform integrability of the
sequence $|g_{k'}(x)-g(x)|^r$.
\endproof

Now we  see that
 $\alpha=0$ and to apply theorem \ref{main_th} we need
 $\chi=\phi+\beta<1$. It is easy to show that the last inequality follows
 from (\ref{subcritfgh}) if only the constant $S$ is sufficiently small and the
 constant $\lambda$ is chosen to make the expression
  $pr$
  to be sufficiently close to $q$.

\subsection{The scale of analytic functions.}

Let
$\mathbb{T}^m=\mathbb{R}^m/(2\pi\mathbb{Z})^m$ be the $m-$dimensional torus.
 All the technique developed
below  can be transferred almost literally to  the case of
the problem with zero boundary conditions on the
$m-$dimensional
cube.

By $x=(x_1,\ldots,x_m)$ denote an element of $\mathbb{R}^m$.

Let
$\mathbb{T}^m_s=\{z=x+iy\in \mathbb{C}^m\mid x\in
\mathbb{T}^m,\quad |y_j|<s,\quad j=1,\ldots,m\}$ be the
complex neighborhood of the torus $\mathbb{T}^m$.


Define a set $E_s,\quad s>0$ as follows
$E_s= C(\overline{\mathbb{T}}^m_s)\bigcap \mathcal{O}(\mathbb{T}^m_s)$.
Here $\mathcal{O}(\mathbb{T}^m_s)$ stands for
the set of analytic functions in  $\mathbb{T}^m_s$.

  The set $E_s$ is a Banach space with respect to the norm
$\|u\|_s=\max_{z\in\overline{\mathbb{T}}^m_s}|u(z)|$. By the Montel
theorem the embeddings $E_{s+\delta}\subset E_s,\quad \delta>0$ are completely continuous.
 By definition put
 $E_0=C(\mathbb{T}^m)$ and $\|\cdot\|_0=\|\cdot\|_{C(\mathbb{T}^m)}.$

Let $\Delta$ stands for the standard Laplace operator
$$\Delta=\sum_{j=1}^m\partial^2_j,\quad
\partial_j=\frac{\partial}{\partial x_j}.$$

\begin{lem}
\label{f345354}
There exists a positive constant $c$ such that for any $u\in E_s,\quad s\ge 0$
the following inequality holds
$$
\|e^{t\Delta}u\|_{s+\delta}\le
c\exp\Big(\frac{\delta^2}{4t}\Big)\|u\|_s,\quad t,\delta>0.$$
The constant $c$ depends only on $m$.
\end{lem}
\proof
The assertion of the lemma easily follows from the well-known
formula:
$$(e^{t\Delta}u)(x)=\frac{1}{(4\pi t)^{m/2}}
\int_{\mathbb{R}}e^{-(\xi_1-x_1)^2/(4t)}\,d\xi_1\ldots
\int_{\mathbb{R}}e^{-(\xi_m-x_m)^2/(4t)}\,d\xi_mu(\xi).$$
In all these integrals one must shift
the  contour of integration to the complex plane and then the desired inequality
follows from the standard estimates.
\endproof

By lemma \ref{f345354} the semigroup $e^{t\Delta}$ is parabolic with
$\gamma=2.$

\begin{lem}
\label{dfg123--}Take a constant $\rho\in (0,1/2].$
For any $\eps\in (0,2\rho)$ there is a positive constant $c=c(\eps)$ such that
if $u\in E_{s+\delta}$ then
\begin{align}
\|(-\Delta)^{-\rho}\partial_j u\|_s&\le
\frac{c}{\delta^{1-2\rho+\eps}}\|u\|_{s+\delta},\quad s\ge
0,\quad\delta>0,\label{vbnvbn}\\
\|(-\Delta)^\rho
u\|_s&\le\frac{c}{\delta^{2\rho+\eps}}\|u\|_{s+\delta}.\label{ety456}\end{align}
\end{lem}
\proof Let us prove formula (\ref{vbnvbn}).
Using the standard facts on Sobolev's spaces we
have
$$\|(-\Delta)^{-\rho}\partial_j u\|_s\le c\|(-\Delta)^{-\rho}\partial_j
u\|_{H^{\eps,p}(\mathbb{T}_s^m)}\le c\|u\|_{H^{\eps+1-2\rho,p}(\mathbb{T}_s^m)},\quad \eps p>2m.$$
Then the desired result follows from the interpolation formula
and
the Cauchy inequality:
$$\|u\|_{H^{\eps+1-2\rho,p}(\mathbb{T}_s^m)}\le
c\|u\|^{\eps+1-2\rho}_{H^{1,p}(\mathbb{T}_s^m)}\|u\|^{2\rho-\eps}_{L^p(\mathbb{T}_s^m)},\quad
\|u\|_{H^{1,p}(\mathbb{T}_s^m)}\le\frac{c}{\delta}\|u\|_{s+\delta}.$$
Formula (\ref{ety456}) is derived in the same way.
\endproof

\begin{prop}[\cite{Taylor}]\label{vtr00}For any constants
 $\quad a\ge r\ge 0$
one has
$$\|e^{t\Delta}u\|_{H^{a}(\mathbb{T}^m)}\le
\frac{c}{t^{(a-r)/2}}\|u\|_{H^{r}(\mathbb{T}^m)}.$$

If  $a>m/2$ then $\|u\|_0\le c \|u\|_{H^{a}(\mathbb{T}^m)}.$
\end{prop}

The first of the following two examples illustrates the effect described in
the Introduction, the second one is to compare our result with the known
one.

\subsection{Integro-differential parabolic equation}In this section we use
the scales $G_s=E_s= C(\overline{\mathbb{T}}^m_s)\bigcap \mathcal{O}(\mathbb{T}^m_s)$.
Let us  focus our
attention on a one dimensional
($m=1$) system.

Consider a problem
\begin{equation}\label{qwqe4575}
u_t=\|(-\Delta)^n u\mid_{\mathbb{T}^m}\|^\lambda_{L^2(\mathbb{T})}+\Delta u,\quad
u\mid_{t=0}=\hat{u}(x)=\sum_{|k|\ge 2}\frac{e^{ikx}}{|k|^{1/2}\log |k|}\in
L^2(\mathbb{T}).\end{equation}
Here $\lambda$ is a positive parameter, $n\in \mathbb{N}$.

Parabolic equations with right side depending on $L^p$ norms of the
unknown function  arise in the theory of
incompressible viscous fluid
\cite{Ohkitani}.

After the change of variable $u=e^{t\Delta}\hat{u}+v$ our problem takes
the form
\begin{equation}\label{zqw5}
v_t=f(t,v)+\Delta v,\quad v\mid_{t=0}=0,\quad
f(t,v)=\|(-\Delta)^n e^{t\Delta}\hat{u}+(-\Delta)^n v\|^\lambda_{L^2(\mathbb{T})}.\end{equation}
One can show that  $f$ is a Lipschitz function:
$$|f(t,v')-f(t,v'')|\le \frac{c}{t^ns^{2n}}\|v'-v''\|_{s}.$$
But this property can not save the situation: the denominator $t^ns^{2n}$
is too small to find a solution by means of a successive procedure. So it
is convenient to  ignore this Lipschitz inequality and write more
effective estimates.

Let us show that if $n\lambda<1$  then problem
(\ref{qwqe4575}) has a solution in the sense of theorem \ref{main_th}.

So that one has
$|f(t,v)|\le
c(\|e^{t\Delta}\hat{u}\|_{H^{2n}(\mathbb{T})}^\lambda+\|(-\Delta)^n v\|^\lambda_{L^2(\mathbb{T})}).$
Then using proposition \ref{vtr00} we obtain
$\|e^{t\Delta}\hat{u}\|_{H^{2n}(\mathbb{T})}\le ct^{-n}\|\hat{u}\|_{L^2(\mathbb{T})} $.
The Cauchy inequality gives
$$\|(-\Delta)^n v\|_{L^2(\mathbb{T})}\le c\|(-\Delta)^n v\|_s\le c\delta^{-2n}\|v\|_{s+\delta},\quad \delta>0. $$
Combining these inequalities with each other
 and taking into account that $(s+\delta)^2<t$
we have
$$|f(t,v)|\le c\delta^{-2n\lambda}(\|\hat{u}\|_{L^2(\mathbb{T})}^\lambda
+\|v\|^\lambda_{s+\delta}).$$ Thus
$\chi=n\lambda$ and if $n\lambda<1$ then by theorem \ref{main_th} the
problem has at least one analytic solution.

Consider the case $\lambda=1$ and let for simplicity $n=1$.

Denote by $u_k$ the Fourier coefficients of a function $u$:
$u(x)=\sum_{k\in \mathbb{Z}}u_ke^{ikx}.$
Notice that the norm of $L^2(\mathbb{T})$ can be presented as
follows
$$\|u\|^2_{L^2(\mathbb{T})}=c\sum_{k\in \mathbb{Z}}|u_k|^2.$$

Then separating the variables in problem (\ref{qwqe4575})
we obtain
\begin{align}
u_0&=c\int_0^t\Big(\sum_{|k|\ge
2}\frac{|k|^3e^{-2\xi|k|^2}}{(\log|k|)^2}\Big)^{\frac{1}{2}}\,d\xi,\label{rty87gh}\\
u_k&=0,\quad \mbox{if}\quad |k|=1,\nonumber\\
u_k&=\frac{e^{-t|k|^2}}{|k|^{1/2}\log|k|},\quad \mbox{if}\quad |k|\ge 2.\nonumber
\end{align}
It is not difficult to show that
$$\Big(\sum_{|k|\ge
2}\frac{|k|^3e^{-2\xi|k|^2}}{(\log|k|)^2}\Big)^{\frac{1}{2}}\ge
-\frac{c}{\xi\log \xi},\quad \xi\in (0,1).$$
So that  the integral in formula (\ref{rty87gh}) does not exist and thus there
are no solutions in this case.

\subsection{3-D Navier-Stokes equation}
In this section we use
the scale $G_s=E_s= C(\overline{\mathbb{T}}^m_s)\bigcap \mathcal{O}(\mathbb{T}^m_s)$.

Consider the Navier-Stokes equation in the divergence free setup.
After  Leray's projection the Navier-Stokes equation takes the well-known form
\begin{equation}
\label{N-S_main}
\begin{aligned}
(u^k)_t&=A^k_l\partial_j(u^ju^l)+\Delta
u^k,\quad A^k_l=(\Delta^{-1}\partial_k\partial_l-\delta_{kl}),\\
u^k\mid_{t=0}&=\hat{u}^k\in H^r(\mathbb{T}^3),
\end{aligned}
\end{equation}
where $\delta_{kl}=1$ for $k=l$ and $0$ otherwise; $k,l,j=1,2,3$
we also use the Einstein
summation convention.

From \cite{kf}, \cite{fk} it follows that if $r=1/2$ then problem
(\ref{N-S_main}) has a solution $u^i(t,x)$ which is regular in the spatial variables for
all $t\in (0,T_*)$. Here  $T_*$ is a small positive constant.

Let us show that by theorem \ref{main_th} the analytic
solution exists for all  $r>1/2$. This indicates that in terms of paper \cite{arrieta}
 theorem \ref{main_th} allows us to carry out only the subcritical case.
This is no surprise since theorem \ref{main_th} is very general.

Assume  a parameter $\rho\in (0,1/2)$ to be close $1/2$
and let us change the variable in (\ref{N-S_main}): $u^k=e^{t\Delta}\hat{u}^k
+(-\Delta)^\rho v^k.$ Then the problem has the form
$$
v^k_t=f^k(t,v)+\Delta
v^k,\quad
v^k\mid_{t=0}=0,$$here\begin{align}
f^k(t,v)&=A_l^k\partial_j(-\Delta)^{-\rho}
(e^{t\Delta}\hat{u}^je^{t\Delta}\hat{u}^l+e^{t\Delta}\hat{u}^j(-\Delta)^\rho
v^l\nonumber\\
&+(-\Delta)^\rho v^je^{t\Delta}\hat{u}^l+(-\Delta)^\rho v^j(-\Delta)^\rho
v^l).\nonumber\end{align}

Estimate the function $f$ term by term. Using lemma \ref{dfg123--} we have
\begin{align}
\|A_l^k\partial_j(-\Delta)^{-\rho}&((-\Delta)^\rho v^j(-\Delta)^\rho v^l)\|_s\le
\frac{c}{\delta^{\eps+1-2\rho}}\sum_{j,l=1}^3\|(-\Delta)^\rho v^j(-\Delta)^\rho
v^l\|_{s+\delta/2}\nonumber\\
&\le \frac{c}{\delta^{\eps+1-2\rho}}\sum_{j,l=1}^3\|(-\Delta)^\rho v^j\|_{s+\delta/2}
\|(-\Delta)^\rho v^l\|_{s+\delta/2}\nonumber\\
&\le\frac{c}{\delta^{\eps+1+2\rho}}\sum_{j,l=1}^3\|v^j\|_{s+\delta}\|v^l\|_{s+\delta}.\nonumber
\end{align}
Now one must choose the parameters $\eps>0, \quad\rho\in (0,1/2)$ such that
\begin{equation}\label{one_____}
\frac{\eps+1+2\rho}{2}<1.\end{equation}

Let us estimate another term of the function $f$ by using lemmas\ref{f345354}, \ref{dfg123--} and
proposition \ref{vtr00} ($(s+\delta)^2<t$):
\begin{align}
\|A_l^k\partial_j(-\Delta)^{-\rho}(e^{t\Delta}\hat{u}^je^{t\Delta}\hat{u}^l)\|_s
&\le\frac{c}{\delta^{1+\eps-2\rho}}\sum_{j,l=1}^3
\|e^{t\Delta}\hat{u}^j\|_{s+\delta}\|e^{t\Delta}\hat{u}^l\|_{s+\delta}\nonumber\\
&\le
\frac{c}{\delta^{1+\eps-2\rho}}\sum_{j,l=1}^3
\|e^{t\Delta/2}\hat{u}^j\|_{0}\|e^{t\Delta/2}\hat{u}^l\|_{0}\nonumber\\
&\le
\frac{c}{\delta^{1+\eps-2\rho}}\sum_{j,l=1}^3
\|e^{t\Delta/2}\hat{u}^j\|_{H^a(\mathbb{T}^3)}\|e^{t\Delta/2}\hat{u}^l\|_{H^a(\mathbb{T}^3)}
\nonumber\\
&\le
\frac{c}{\delta^{1+\eps-2\rho}t^{a-r}}\sum_{j,l=1}^3
\|\hat{u}^j\|_{H^r(\mathbb{T}^3)}\|\hat{u}^l\|_{H^r(\mathbb{T}^3)},
\nonumber
\end{align}here $a>3/2$. We need to have
\begin{equation}
\label{two________}
\frac{1+\eps-2\rho}{2}+a-r<1.\end{equation}
In the same manner we obtain
$$\|A_l^k\partial_j(-\Delta)^{-\rho}(e^{t\Delta}\hat{u}^j(-\Delta)^\rho v^l)\|_s
\le \frac{c}{\delta^{\eps+1}t^{(a-r)/2}}
\sum_{j,l=1}^3\|\hat{u}^j\|_{H^r(\mathbb{T}^3)}\|v^l\|_{s+\delta}.$$
Thus there must be
\begin{equation}
\label{fhgr7869}
\eps+1+a-r<2.\end{equation}
It is not difficult to show that for any $r>1/2$ there exists the small
parameter $\eps>0$, the  parameter $a$ close to $3/2$ from above and the parameter
$\rho$ close to $1/2$ from below
such that inequalities (\ref{one_____}), (\ref{two________}),
(\ref{fhgr7869}) are fulfilled.


 \end{document}